\documentclass[11pt]{amsart}

\usepackage[all]{xypic}

\usepackage{latexsym}

\usepackage{amssymb}

\usepackage{amsfonts}

\usepackage{amscd}

\usepackage{amsmath,amsthm}

\newtheorem{lemma}{Lemma}[section]

\newtheorem{theorem}[lemma]{Theorem}

\newtheorem{lem}[lemma]{Lemma}

\newtheorem{prop}[lemma]{Proposition}

\newtheorem{thm}[lemma]{Theorem}

{

}

\theoremstyle{definition}

\theoremstyle{remark}

\numberwithin{equation}{section}

\newenvironment{pf}{\noindent{\bf Proof.}}{\hfill $\square$\medskip}

\def\NN{{\mathbb N}}

\def\PP{{\mathbb P}}

\def\ZZ{{\mathbb Z}}

\def\fol{{\bar f}}

\def\Iol{{\bar I}}

\def\Sol{{\bar S}}

\def\0ol{{\bar 0}}

\def\1ol{{\bar 1}}

\def\2ol{{\bar 2}}

\def\ol2{{\bar 2}}

\def\3ol{{\bar 3}}

\def\4ol{{\bar 4}}

\def\5ol{{\bar 5}}

\def\6ol{{\bar 6}}

\def\7ol{{\bar 7}}

\def\8ol{{\bar 8}}

\def\9ol{{\bar 9}}

\def\bold0{{\bf 0}}

\def\bold1{{\bf 1}}

\def\bold2{{\bf 2}} 

\def\bold3{{\bf  3}}

\def\bold4{{\bf 4}}

\def\bold5{{\bf 5}}

\def\bold6{{\bf 6}}

\def\bold7{{\bf 7}}

\def\bold8{{\bf 8}}

\def\bold9{{\bf 9}}

\def\P2Skly{\PP^2_{Skly}}

\def\coker{\operatorname {coker}}

\def\ker{\operatorname {ker}}

\def\th{\operatorname {th}}    

\def\Fdim{{\sf Fdim}}

\def\Gr{{\sf Gr}}

\def\max{\operatorname{max}}

\def\QGr{\operatorname{\sf QGr}}

\def\ul1{\operatorname{\underline{1}}}

\def\dirlim{\mathop{\vtop{\baselineskip -100pt\lineskip -1pt\lineskiplimit 0pt

\setbox0\hbox{lim}\copy0\hbox to \wd0{\rightarrowfill}}}\limits}

\def\invlim{\mathop{\vtop{\baselineskip -100pt\lineskip -1pt\lineskiplimit 0pt

\setbox0\hbox{lim}\copy0\hbox to \wd0{\leftarrowfill}}}\limits}

\def\I11{{1 \kern -0.8pt \! \mbox{l}}}

\def\mumu{{\mu\kern-4.2pt\mu}}

\def\bfmu{{\mu\kern-4.2pt\mu}}

\def\2slash{\backslash \! \backslash}

\def\boxtimes{\setbox0\hbox{$\Box$}\copy0\kern-\wd0\hbox{$\times$}}

\date{}                                           

\begin{document}

\title[Graded modules over monomial algebras and path algebras]{An equivalence of categories for 
graded modules over monomial algebras and path algebras of quivers}

\author{Cody Holdaway and S. Paul Smith}

\address{ Department of Mathematics, Box 354350, Univ.
Washington, Seattle, WA 98195}

\email{codyh3@math.washington.edu, smith@math.washington.edu}




\begin{abstract}

Let $A$ be a finitely generated connected graded $k$-algebra defined  by a finite number of 
monomial relations, or, more generally, the path algebra of a finite quiver modulo a finite number of relations of the form ``path$=0$''. Then there is a finite directed graph, $Q$, the Ufnarovskii graph of $A$, for 
which there is an equivalence of categories $\QGr A \equiv \QGr(kQ)$. Here $\QGr A$ is the quotient category  
$\Gr A/\Fdim$ of graded $A$-modules modulo the subcategory consisting of those that are the sum of their finite dimensional submodules. The proof makes use of an algebra homomorphism $A \to kQ$ that may be of independent interest. 
\end{abstract}

\maketitle

\pagenumbering{arabic}


\setcounter{section}{0}

\section{Introduction}

\subsection{}

Throughout $k$ is a field.

Let $A$ be an $\NN$-graded $k$-algebra.

The category of $\ZZ$-graded {\it right} $A$-modules with degree-preserving homomorphisms
is denoted by $\Gr A$ and  $\Fdim A$ is its full subcategory consisting of modules that are the sum of their
finite-dimensional submodules. 
Since $\Fdim A$ is a Serre subcategory of $\Gr A$ (it is, in fact, a localizing subcategory) we may form the 
quotient category 
$$
\QGr A:=\frac{\Gr A}{\Fdim  A}.
$$

We are interested in the structure of $\QGr A$ for monomial algebras.

\subsection{}
A connected graded {\sf monomial algebra} is a free algebra modulo an ideal generated by words in the letters generating the free algebra. More explicitly, if $w_1,\ldots,w_r$ are words in the letters $x_1,\ldots,x_g$, then
\begin{equation}
\label{A.defn.1}
A=\frac{k \langle x_1,\ldots,x_g\rangle}{(w_1,\ldots, w_r)}
\end{equation}
is a  finitely presented  monomial algebra.

Our main result applies to a more general class of monomial algebras, namely those of the form $kQ'/I$
where $Q'$ is a finite quiver (section \ref{sect.notn}) and $I$ an ideal generated by a finite set of paths in $Q'$. 
Such algebras can be described without mentioning quivers: let $K$ be a finite product of copies of $k$, $T_KV$ the tensor algebra of a $K$-bimodule $V$ that has a finite $k$-basis $a_1,\ldots,a_g$, and
\begin{equation}
\label{A.defn.2}
A=\frac{T_KV}{(p_1,\ldots,p_r)}
\end{equation}
where each $p_j$ is a word in the $a_i$'s.

\subsection{The main result}

\begin{thm}
\label{thm.main0}
Let $A$ be a monomial algebra of the form (\ref{A.defn.2}).  There is a quiver $Q$ and 
an equivalence of categories 
$$
\QGr A \equiv \QGr kQ.
$$
\end{thm}

The structure and properties of  $\QGr kQ$ are  described in \cite{Sm2}.

The proof of Theorem \ref{thm.main0} uses result of Artin and Zhang, Proposition \ref{prop.AZ}
below, in an essential way.

When $A$ is of the form (\ref{A.defn.1}) we can take $Q$ to be its Ufnarovskii graph (section \ref{sect.Ufn})
and there is then a  homomorphism $f:A \to kQ$ such that the functor $- \otimes_A kQ$ induces the equivalence in Theorem \ref{thm.main0}. This is proved in section \ref{sect.4.1}; see Theorem \ref{thm.main} for a precise statement.

In section \ref{sect.4.2}, Theorem \ref{thm.main0} is proved for algebras of the form  (\ref{A.defn.2}): 
if $A$ is of the form  (\ref{A.defn.2}) its subalgebra generated by $k$ and $A_1$ is of the form (\ref{A.defn.1}) and has finite codimension in $A$ so, by Artin and Zhang's result and Theorem \ref{thm.main0} for algebras of the form (\ref{A.defn.1}), Theorem \ref{thm.main0} holds for algebras of the form  (\ref{A.defn.2}).

\subsection{Quadratic monomial algebras}

If $A$ is  monomial algebra of the form  (\ref{A.defn.1}) with $\deg w_i=2$ for all $i$ we call $A$ a {\it quadratic}
monomial algebra.  The proof of Theorem \ref{thm.main0} for quadratic monomial algebras is much simpler than the general case. We give that proof in section \ref{sect.quad.algs.1}. 

Let $A$ be an arbitrary finitely presented connected graded monomial algebra. By Backelin-Fr\"oberg \cite{BF}, the Veronese subalgebra $A^{(n)}\subset A$ is quadratic for $n \gg 0$; by Verevkin \cite{V}, $\QGr A \equiv \QGr A^{(n)}$, so Theorem \ref{thm.main0} holds for $A$ if it holds for $A^{(n)}$. However, if \ref{thm.main0}
is proved for $A$ by first proving it for $A^{(n)}$ the quiver $Q$ is the Ufnarovskii graph for $A^{(n)}$ 
which is more complicated than that for $A$ (see section \ref{ssect.example} for an example illustrating  this).  

That is why we prove theorem \ref{thm.main0} directly in section \ref{sect.4.1}, i.e., without passing to a Veronese subalgebra.

\subsection{Acknowledgements}

We thank Victor Ufnarovskii for a comment that prompted us to change our notation and thereby
simplify the proofs. We also thank Chelsea Walton for suggesting several improvements to an earlier 
version of this paper.

\section{Preliminaries}

\subsection{Notation}
\label{sect.notn}

\label{sect.1.1}

The letter $Q$ will always denote a directed graph, or quiver, with a finite number of vertices and arrows---loops and multiple arrows between vertices are allowed.

We write $kQ$ for the path algebra of $Q$. The finite paths in $Q$, including the trivial paths at each vertex,
form a basis for $kQ$ and multiplication is given by concatenation of paths. If $a$ is an arrow that ends where the arrow $b$ begins we write
$$
ab:=\hbox{the path ``$a$ followed by $b$''}.
$$
We set $ab=0$ if $b$ does not begin where $a$ ends.
Likewise, if a path $p$  ends where a path $q$ begins, $pq$ denotes the path {\it first traverse $p$ then $q$}.

We make $kQ$ an $\NN$-graded algebra by declaring that a path is homogeneous of degree equal to its length.  


\subsection{}

Throughout, modules are {\it right} modules.

\begin{prop}
\cite[Prop 2.5]{AZ}
\label{prop.AZ}
Let $\phi:A \to B$ be a homomorphism of graded $k$-algebras. If $\ker \phi$ and $\coker\phi$ belong to $\Fdim A$, then $- \otimes_A B$ induces an equivalence of categories $\QGr A \to \QGr B$.
\end{prop}

\begin{lemma}
\label{lemma1.2}
Let $A$ and $B$ be $\NN$-graded $k$-algebras generated by $A_0+A_1$ and $B_0+B_1$ respectively. 
Let $\phi:A\to B$ be a homomorphism of graded $k$-algebras. If $B_0\phi(A_{m})\subset \phi( A_m)$ and 
$B_1 \phi(A_m)\subset \phi(A_{m+1})$ for some $m \in \NN$, then $\coker \phi$ belongs to  $\Fdim A$.
\end{lemma}
\begin{pf}
We can replace $A$ by its image in $B$ so we will do that; i.e., without loss of generality, $A$ is a graded subalgebra of $B$ and $\phi$ is the inclusion map.

If $n \ge 2$ and $B_{n-1}A_{m}\subset A_{m+n-1}$, then
$$
B_nA_{m}=B_1B_{n-1}A_{m}\subset B_1A_{m+n-1}=B_1A_mA_{n-1}\subset A_{m+1}A_{n-1}=A_{m+n}.
$$
It follows that $B_nA_{m}\subset A_{m+n}$ for all $n\geq 0$.  Thus $B/A$ is annihilated on the right by 
$A_m$ and therefore belongs to  $\Fdim A$.
\end{pf}


\section{The Ufnarovskii graph of a connected graded monomial algebra}
\label{sect.Ufn}

Throughout this paper $G$ is a fixed finite set of {\it letters} or {\it generators}, $\langle G \rangle$ is the free monoid generated by $G$, and $k\langle G\rangle$ is the free $k$-algebra generated by $G$.
Elements of $\langle G \rangle$ are called {\sf words}.
Throughout, $F$ denotes a fixed {\it finite} set of words and 
\begin{equation}
\label{defn.A.notn}
A:=\frac{k\langle G \rangle}{(F)}
\end{equation}
is the quotient by the ideal $(F)$ generated by $F$. Such $A$ is called 
a {\sf monomial algebra}.

There is no loss of generality in assuming that $G \cap F=\varnothing$. We will make that assumption.

We make $A$ a graded algebra by placing $G$ in degree one. Thus $A_1=kG$.

\subsection{Words} 
\label{ssect.words}

The words in $F$ are said to be {\sf forbidden}. 
A word is {\sf illegal} if it belongs to $(F)$ and {\sf legal} otherwise. The set of legal words is denoted by $L$,
and $L_r:=L \cap G^r$ is the set of legal words of length $r$. The image of $L_r$ in $A$ is a basis for $A_r$; see, for example, \cite[Lemma 2.2]{KL}. 

Throughout we use the notation
\begin{align*}
\ell+1:= & \hbox{the longest length of a forbidden word} 
\\
= & \max\{ \ell+1 \; | \; F \cap G^{\ell +1}\ne \varnothing\}, \qquad \hbox{and}
\\
L_{\le r}:= & \{ \hbox{legal words of length $\le r$}\}.
\end{align*}

\subsection{Notation}
The letters $s$, $t$, $u$, $v$, $w$, will always denote words.

If $u$ and $w$ are words we write $$u \lhd w$$ if $w=uv$ for some word $v$. 

The symbols $x$, $y$, and $x_i$, will always denote elements of $G$. The notation $x_i \lhd w$
therefore means that $x_i$ is the first letter of $w$. 

\subsection{The Ufnarovskii graph}
\label{ssect.ufn}
The {\sf Ufnarovskii graph} of $A$ is the directed graph $Q$, or $Q(A)$ if we need to specify $A$, defined as follows (see  \cite[Sect. 12.2]{KL}, \cite{U1},  \cite{U2}). 

The set of vertices of $Q$ is
$$
Q_0=L_\ell.
$$
The set of arrows of $Q$ is in bijection with the set $L_{\ell+1}$ as follows,
$$
Q_1=\{a_w\;|\; w\in L_{\ell+1}\}.
$$ 
If $w \in L_{\ell+1}$, then there are unique $s, t \in Q_0$ and unique $x,y \in G$ such that 
$w=sy = xt  \in L$ and we declare that the arrow $a_w$ corresponding to $w$ goes from $s$ to $t$.  

Given $s,t \in Q_0$, there is at most one arrow from $s$ to $t$.

Suppose  $n>0$. If $x_1\ldots x_{n+\ell}$ is a legal word of length $n+\ell$ there is a length-$n$ path  
\begin{equation}
\label{a.path}
x_1\ldots x_\ell \, \longrightarrow  \, x_2\ldots x_{\ell+1} \, \longrightarrow\, \cdots  \, \longrightarrow \,
x_{n+1} \ldots x_{n+\ell}
\end{equation}
in $Q$. 
This provides a bijection between legal words of length 
$n+\ell$ and paths of length $n$ (see the proof of \cite[Thm. 3]{U1} and the remark at \cite[p.157]{KL}).

 \subsection{Labeling arrows and paths}
 \label{sect.labels}
 
 We write $a_w$ for the arrow  corresponding to $w \in L_{\ell+1}$. The path in (\ref{a.path}) is therefore
$a_{x_1\ldots x_{\ell+1}} a_{x_2\ldots x_{\ell+2}} \cdots a_{x_n\ldots x_{n+\ell}}$.
 
 Suppose there is an arrow $u \to v$.  Then $uy=xv$ for unique $x$ and $y$ in $G$, and we attach the {\sf label}
 $x$ to the arrow $u \to v$. We denote this by $u \stackrel{x}{\longrightarrow} v$. The following facts will be used often:
 \begin{itemize}
  \item 
The label attached to the arrow $a_w$ is the first letter of $w$.
  \item 
  The existence of an arrow $u \stackrel{x}{\longrightarrow} v$ implies that $x \lhd u$ and $u\lhd xv$.  
\end{itemize}

We extend the labeling to paths: the {\sf label} attached to a concatenation of arrows is the concatenation of the labels attached to the arrows in the path---for example, the label attached to the path in (\ref{a.path}) is $x_1\ldots x_n$. In general, there will be different paths with the same label: for example, the labels on the Ufnarovskii graph for $A=k\langle x,y\rangle/(y^3)$ are
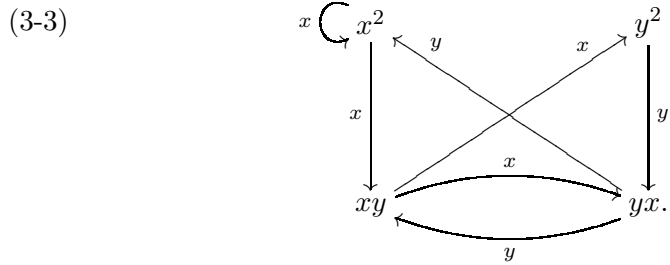
\begin{equation}
\label{quiv.Q}
  \UseComputerModernTips
\xymatrix{ 
\ar@(ul,dl)[]_{x}   x^2    \ar[dd]_{x} &&& y^2  \ar[dd]^{y} 
\\
\\
xy   \ar@/^1pc/[rrr]^{x} \ar[uurrr]^(0.8){x}   &&& yx.  \ar[uulll]_(0.8){y} \ar@/^1pc/[lll]^{y}
}
\end{equation}

The Ufnarovskii graph for  $k\langle x,y,z \rangle/(z^2,zy)$ appears in Section \ref{sect.anexample}.

The following observation is surely known to the experts.

\begin{lemma}
\label{lem.paths}
  Suppose there is a path with the label $x_1\ldots x_r$, say
\begin{equation}
\label{a.path.2}
v_0 \stackrel{x_1}{\longrightarrow} v_1 \stackrel{x_2}{\longrightarrow} \cdots  \stackrel{x_r}{\longrightarrow} v_r.
\end{equation}
Let $v_r=x_{r+1} \ldots x_{r+\ell}$. 
\begin{enumerate}
  \item 
 $v_{i-1} =  x_i\ldots x_{i+\ell-1}$ for all $i=1,\ldots,r+1$. 
 \item{}
$x_1\ldots x_rv_r$ is a legal word.
 \item
  $x_1\ldots x_r \notin (F)$.
\end{enumerate}
\end{lemma}
\begin{pf}
The hypothesis implies $v_{i-1} \lhd x_iv_i$ and $x_i \lhd v_{i-1}$ for all $i=1,\ldots,r$.
An induction argument, or simply noticing the pattern in the equalities
\begin{align*}
v_r & = \, x_{r+1} \ldots x_{r+\ell},
\\
v_{r-1} & = \, x_{r} \ldots x_{r+\ell-1},
\\
v_{r-2} & = \, x_{r-1} \ldots x_{r+\ell-2},  \qquad \hbox{etc.}
\end{align*}
proves (1).

(2)
To prove $x_1\ldots x_rv_r$ is legal it suffices to show its subwords of length $\ell+1$ are legal. 
Such a subword is of the form $ x_{i} \ldots x_{i+\ell-1}x_{i+\ell}$ for some $i$ in the range $1 \le i \le r$; this subword is equal to $v_{i-1}x_{i+\ell} = x_i v_i$ and is legal because there is an arrow $v_{i-1} \to v_i$.

 (3)
Since a subword of a legal word is legal, (3) follows from (2).
\end{pf}

The contrapositive of part (3) of Lemma \ref{lem.paths} is  useful so we record it separately.

\begin{lemma}
\label{lem.paths.2'}
If $x_1\ldots x_r$ is an illegal word, then there are no paths labeled $x_1\ldots x_r$.
\end{lemma}
 
The converse of Lemma \ref{lem.paths.2'} is false. For example, 
$x$ is a legal word when $A=k\langle x,y \rangle/(xy,x^2)$ but the Ufnarovskii graph of $A$ is 
$$
Q=\xymatrix{ {y} \ar@(ul,dl)[]_{a_{yy}} \ar[r]^{a_{yx}} & {x}}
$$
 with labels 
\begin{equation}
\label{y^2=0}
\xymatrix{ {y} \ar@(ul,dl)[]_{y} \ar[r]^{y} & {x.}}
\end{equation}
\vskip .2in

\subsection{The homomorphism $k\langle G\rangle/(F) \to kQ$}

Let $f:k\langle G\rangle \to kQ$ be the unique algebra homomorphism 
such that for all $x \in G$,
$$
f(x)=\hbox{the sum of all arrows labeled $x$}
$$
or $0$ if there are no arrows labeled $x$. 

Hence, if $x_1\ldots x_r \in G^r$, 
\begin{equation}
\label{eq.f(w)}
f(x_1\ldots x_r)=\hbox{the sum of all paths labeled $x_1\ldots x_r$}
\end{equation}
or $0$ if there are no such paths.
More formally,
\begin{align*}
f(x)=&\phantom{xi}  0  \phantom{xxxxxx} \hbox{if $xL_\ell \cap L_{\ell+1} = \varnothing$, and}
\\
f(x)=	& \sum_{\substack{w \in Q_1  \\ x \lhd w}} a_{w}   \phantom{xx} \hbox{if $xL_\ell \cap L_{\ell+1}  \ne \varnothing$.}
\end{align*}
Since  $f(G) \subset Q_1$, $f$  is a homomorphism of graded $k$-algebras.

\begin{prop}
 \label{inducedmap}
The homomorphism $f:k\langle G\rangle \to kQ$ induces a homomorphism of graded algebras 
from $A$ to $kQ$.
\end{prop}
\begin{pf} 
Lemma \ref{lem.paths.2'} and (\ref{eq.f(w)}) imply $f(w)=0$ for all $w \in F$. 
\end{pf}

\begin{lem}
\label{lem.paths.2}
Let $x_1\ldots x_r \in G^r$. There is a path labeled $x_1\ldots x_r$ if and only if $x_1\ldots x_r L_\ell \cap L \ne \varnothing$.
\end{lem}
\begin{pf}
($\Rightarrow$)
Suppose there is a path
$$
v_0 \stackrel{x_1}{\longrightarrow} v_1 \stackrel{x_2}{\longrightarrow} \cdots  \stackrel{x_r}{\longrightarrow} v_r.
$$
Write $v_r=x_{r+1}\ldots x_{r+\ell}$. Since $x_iv_i$ is legal for all $i=1,\ldots,r$ and $x_iv_i=x_ix_{i+1} \ldots x_{i+\ell -1}$ all subwords of $x_1\ldots x_r v_r$ of length $\ell+1$ are legal. It follows that $x_1\ldots x_r v_r$ is legal.

($\Leftarrow$) 
Suppose $x_1\ldots x_r L_\ell  \cap L \ne \varnothing$. Let $v_r=x_{r+1} \ldots x_{r+\ell}$ be a vertex such that
$x_1\ldots x_rv_r$ is legal. For $i=1,\ldots,r$, define
$$
v_{i-1}:= x_i \ldots x_{i+\ell -1}.
$$
This is a legal word, of length $\ell$, because it is a subword of the legal word $x_1\ldots x_rv_r$. Since 
$v_{i-1} \lhd x_i v_i$ there is an arrow $v_{i-1} \stackrel{x_i}{\longrightarrow} v_i$. 
Concatenating these arrows produces a 
path labeled $x_1\ldots x_r$.
\end{pf}

\begin{lem}
\label{lem.paths.3}
Let $x_1\ldots x_r$ be a legal word of length $r  \ge \ell$. 
There is a path labeled $x_1\ldots x_r$ if and only if there is a path labeled $x_{r-\ell+1} \ldots x_r$.
\end{lem}
\begin{pf}
The lemma is true for $r=\ell$ so suppose $r>\ell$. 

($\Rightarrow$)
This is obvious.

($\Leftarrow$) 
Suppose there is a path
$$
v_{r-\ell} \stackrel{x_{r-\ell+1}}{\longrightarrow} v_{r-\ell +1} \longrightarrow \cdots  \longrightarrow
v_{r-1} \stackrel{x_r}{\longrightarrow} v_r.
$$
Write $v_r=x_{r+1}\ldots x_{r+\ell}$. 

By Lemma \ref{lem.paths.2}, $x_1\ldots x_r$ is legal if $x_1\ldots x_rv_r$ is. The word $x_1\ldots x_rv_r$
is legal if all its 
subwords of length $\ell+1$ are legal. The proof of Lemma \ref{lem.paths.2}
showed that $x_{r-\ell+1} \ldots x_rv_r$ is legal. All subwords of $x_{r-\ell+1} \ldots x_rx_{r+1}\ldots x_{r+\ell}$
are therefore legal so it only remains to show that $x_i \ldots x_{i+\ell}$ is legal for all $i \le r-\ell$. If 
 $i \le r-\ell$, then $x_i \ldots x_{i+\ell}$ is a subword of $x_1\ldots x_r$ and therefore legal.  
\end{pf}

\subsection{The kernel of $f$}

The homomorphism $f$ need not be injective: for example, by looking at the labels on the 
quiver (\ref{y^2=0}) above one sees that $f(x)=0$ when $A=k\langle x,y \rangle/(xy,x^2)$.  

\begin{lemma}
 \label{injonwords}
Let $w_1,\ldots,w_n$ be pairwise distinct legal words. If $f(w_i) \ne 0$ for all $i$, then 
$\{f(w_1),\ldots,f(w_n)\}$ is linearly independent.
\end{lemma}
\begin{pf}
Since $f$ preserves degree we can assume that $w_1,\ldots,w_n$ have the same length, say $r$.  By definition, $f(w_i)$ is the sum of the paths labeled $w_i$; hence if $i \ne j$ no path that appears in $f(w_i)$ appears in $f(w_j)$. But the paths of length $r$ are linearly independent elements of $kQ$ so $\{f(w_1),\ldots,f(w_n)\}$ is linearly independent.
\end{pf}

\begin{thm}
\label{thm.ker}
The kernel of the homomorphism $f:k\langle G \rangle \to kQ$ is equal to $(F)+I$ where $I$ is 
the left ideal generated by the set
$$
S:=\{x_1\ldots x_s \in G^s \; | \; s \le \ell \; \hbox{and there is no path labeled $x_1\ldots x_s$}\}.
$$ 
\end{thm}
\begin{pf}
By Proposition \ref{inducedmap}, $\ker f$ contains the ideal $(F)$.
Since $f(x_1\ldots x_r)$ is the sum of all the paths labeled $x_1\ldots x_r$, $S \subset \ker f$.
Hence $(F)+I \subset \ker f$. 

Since $(F)$ is spanned by words, Lemma \ref{injonwords} implies $\ker f$ is spanned by $(F)$ and various
legal words. 
Suppose $x_1\ldots x_r$ is a legal word such that $f(x_1\ldots x_r)=0$. 
This implies there is no path labeled $x_1\ldots x_r$ so, if $r \le \ell$, $x_1\ldots x_r$ is in $S$ and therefore in $I$. 
On the other hand, if $r \ge \ell+1$, Lemma \ref{lem.paths.3} implies $x_{r-\ell+1} \ldots x_r$ is in $S$, whence 
$x_{1} \ldots x_r \in I$.  
\end{pf}

Information about the cokernel of $f$ is given in Proposition \ref{kerinfdim}.


\section{The Proof of Theorem 1.1}  
\label{sect.thmmain}

\subsection{The proof of Theorem \ref{thm.main0} when $A$ is as in (\ref{A.defn.1})}
\label{sect.4.1}

Let $A$ be as in (\ref{A.defn.1}) and adopt the notation in (\ref{defn.A.notn}). 
We will prove Theorem \ref{thm.main0} by applying Proposition \ref{prop.AZ} to the induced homomorphism 
$\fol:A \to kQ$.  Before doing that we must check that the hypotheses of Proposition \ref{prop.AZ}  hold:
 we must show that the kernel  and cokernel of $\fol$ belong to $\Fdim A$.

\begin{prop} \label{kerinfdim}
Let $\fol:A\to kQ$ be the homomorphism induced by $f$. Then $\ker \fol$ and $\coker \fol$ belong to $\Fdim A$.
\end{prop}
\begin{pf}
Let $I$ and $S$ be as in Theorem \ref{thm.ker} and write $\Iol$ and $\Sol$ 
for their images in $A$. Thus, $\Iol=\ker \fol$ and $\ker \fol$ is generated as a left ideal by $\Sol$. 

Given the description of $\ker f$ in Theorem \ref{thm.ker}, it suffices to show that $\Iol A_{\ell}=0$. 

Let $x_1\cdots x_s \in S$. By Lemma \ref{lem.paths.2}, $x_1\ldots x_r L_\ell \cap L = \varnothing$;
in other words, $x_1\ldots x_r L_\ell \subset (F)$. Taking the image of this equality in $A$ we conclude that $\Sol A_{\ell}=0$. It follows that $\Iol A_{\ell}=0$. Thus $\ker f$ belongs to $\Fdim A$.

By Lemma \ref{lemma1.2}, to show $\coker \fol$ belongs to $\Fdim A$ it suffices to show that
$$
(kQ_0) \fol(A_\ell) \subset  \fol(A_\ell)
\qquad \hbox{and} \qquad
(kQ_1) \fol(A_\ell) \subset  \fol(A_{\ell+1}).
$$
To do this it suffices to show that $Q_0f(L_\ell) \subset f(L_\ell)$ and $Q_1f(L_\ell) \subset f(L_{\ell+1})$.

Let $x_1\ldots x_\ell \in L_\ell$. 
By Lemma \ref{lem.paths}(1), every path labeled $x_1\ldots x_\ell$ begins
at the vertex $v_0= x_1\ldots x_\ell$. 

Let $e$ be a trivial path and $p$ a path labeled $x_1\ldots x_\ell$; 
since $p$ begins at $v_0$, $ep=p$ if $e$ is the trivial path  at $v_0$ ,and $ep=0$ if $e$ is some other trivial path.
Hence $ef(x_1\ldots x_\ell)$ is either 0 or  $f(x_1\ldots x_\ell)$. It follows that $Q_0f(x_1\ldots
x_{\ell}) = \{f(x_1\ldots x_{\ell})\}$ and $Q_0f(L_\ell) = f(L_\ell)$.

Let $a$ be an arrow and $p$ a path labeled $x_1\ldots x_\ell$. If $a$ does not end at $v_0$, then $ap=0$ because $p$ begins at $v_0$; thus, if  $a$ does not end at $v_0$, then $af(x_1\ldots x_\ell)=0$.

We now assume $a$ ends at $v_0$; i.e., $v_{-1} \stackrel{a}{\longrightarrow} v_0$ and the arrow $a$ is labeled by the first letter of $v_{-1}$, say $x_{0}$. The path $ap$ is therefore labeled $x_0x_1 \ldots x_\ell$. 
Since $v_0\lhd x_0v_1$, $a$ is the only arrow labeled $x_0$ that ends at $v_0$. Therefore
\begin{align*}
af(x_1\ldots x_\ell) = & f(x_0)f(x_1\ldots x_\ell)
\\
= & f(x_0x_1\ldots x_\ell).
\end{align*}
In particular, $af(x_1\ldots x_\ell)  \in f(L_{\ell+1})$. 

This completes the proof that $Q_1f(L_{\ell}) \subset f(L_{\ell+1})$ and, as explained before, this implies $\coker
\fol$ belongs to $\Fdim A$.
\end{pf}

\begin{theorem}
\label{thm.main}
Let $A$ be a connected graded monomial algebra as in (\ref{A.defn.1}) and/or (\ref{defn.A.notn}).
Let $Q$ be its Ufnarovskii graph and view $kQ$ as a left $A$-module through the homomorphism 
$\fol:A \to kQ$. Then $-\otimes_A kQ$ induces an equivalence of categories $\QGr A\equiv \QGr kQ$.
\end{theorem}
\begin{pf}
This follows from Propositions  \ref{prop.AZ} and \ref{kerinfdim}.
\end{pf}

\subsection{The proof of Theorem \ref{thm.main0} when $A$ is as in (\ref{A.defn.2})}
\label{sect.4.2}

Let $Q'$ be a finite quiver and $A=kQ'/I$ the quotient of its path algebra by an ideal generated by a finite number of paths. (Thus $A$ is a more general kind of monomial algebra.) The subalgebra 
$$
A'=k \oplus A_1 \oplus A_2 \oplus \cdots
$$
is of finite codimension in $A$ so $A/A' \in \Fdim A'$. Proposition  \ref{prop.AZ} therefore implies that
$- \otimes_{A'}A$ induces an equivalence of categories
\begin{equation}
\label{equiv1}
\QGr A' \equiv \QGr A.
\end{equation}

Since $A'$ is a monomial algebra of the form (\ref{A.defn.1}),  Theorem \ref{thm.main} gives an equivalence
\begin{equation}
\label{equiv2}
\QGr A' \equiv \QGr kQ
\end{equation}
where $Q$ is the Ufnarovskii graph of $A'$. By (\ref{equiv1}) and (\ref{equiv2}),
$$
\QGr A \equiv \QGr kQ.
$$

This completes the proof of Theorem \ref{thm.main0} for $kQ'/I$.


\section{An Example}
\label{sect.anexample}

Let $A=k\langle x,y,z\rangle/(z^2,zy)$. Since $\ell=1$, $Q_0=\{x,y,z\}$. The arrows for $Q(A)$ correspond to the legal words of length two, namely
$$
\{x^2,xy,xz,y^2,yx,yz,z^2,zx,zy\}-\{z^2,zy\}.
$$
The Ufnarovskii graph of $A$ is therefore
$$
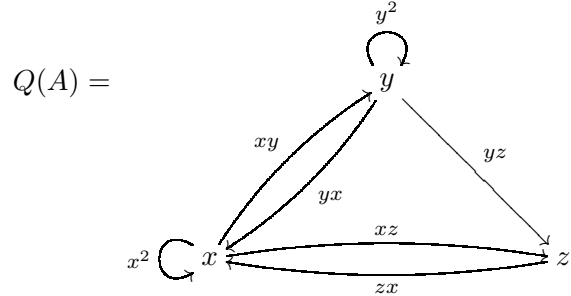
\UseComputerModernTips
Q(A)=\xymatrix{ & & {y} \ar@(ul,ur)[]^{y^2} \ar[ddrr]^{yz} \ar@/^/[ddll]^{yx} & & \\
 & & & & \\
          {x} \ar@(ul,dl)[]_{x^2} \ar@/^/[rrrr]^{xz} \ar@/^/[rruu]^{xy} & & & & {z} \ar@/^/[llll]^{zx} }
$$
(the arrows are denoted by $w$ rather than $a_w$) with labels
$$
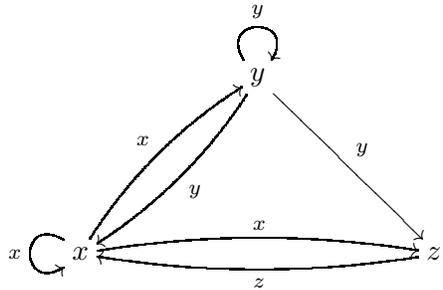
\UseComputerModernTips
\xymatrix{ & & {y} \ar@(ul,ur)[]^{y} \ar[ddrr]^{y} \ar@/^/[ddll]^{y} & & \\
 & & & & \\
          {x} \ar@(ul,dl)[]_{x} \ar@/^/[rrrr]^{x} \ar@/^/[rruu]^{x} & & & & {z.} \ar@/^/[llll]^{z} }
$$
Thus, the homomorphism $f$ is 
\begin{eqnarray*}
f(x)&=&a_{x^2}+a_{xy}+a_{xz} \\
f(y)&=&a_{y^2}+a_{yx}+a_{yz} \\
f(z)&=&a_{zx}.
\end{eqnarray*}


\section{Connected graded quadratic monomial algebras}
\label{sect.quad.algs}

Section \ref{sect.quad.algs.1} contains a short proof of Theorem \ref{thm.main} for connected graded monomial algebras with quadratic relations. Section \ref{sect.quad.algs.2} shows that Theorem \ref{thm.main}
for an arbitrary finitely presented connected graded monomial algebra $A$ can be deduced from the 
quadratic case.

\subsection{}
\label{sect.quad.algs.1}
Let $A$ be a quadratic monomial algebra and $Q$ its Ufnarovski graph. 

The defining relations for $A$ have length $2$ so $\ell=1$. The set of vertices for $Q$ is therefore in bijection with $G$. There is an arrow $a_{xy}$ from vertex $x$ to vertex $y$ if and only if $xy \notin F$ and that arrow is labeled  
$x$ if it exists. It follows that the map $f:k\langle G\rangle \to kQ$ defined in section 3 can be defined as follows:  
$$
f(x)= \hbox{the sum of all arrows that start at $x$}.
$$  
Thus, if $r \ge 2$, then 
$$
f(x_1\ldots x_r)=
	\begin{cases}
		pf(x_r) & \text{where $p$ is the unique path labeled}
		\\
		& \hbox{ $x_1\ldots x_{r-1}$ that ends at vertex $x_r$;}
		\\
		0  & \text{if there is no such $p$.}
	\end{cases}	
$$		
In particular, if $xy \in F$, there is no arrow from $x$ to $y$ so $f(xy)=0$. Thus $f(F)=0$ and there is an induced 
map $\fol:A \to kQ$. 

The lemmas in section \ref{sect.Ufn} are either trivial or unnecessary in the quadratic case. 
The proof that $\ker \fol$ belongs to $\Fdim A$ is also much simpler.

\subsection{}
\label{sect.quad.algs.2}
Let $n$ be a positive integer. The $n^{\th}$ {\sf Veronese subalgebra} of a $\ZZ$-graded algebra $B$ is
$$
B^{(n)}:= \bigoplus_{i \in \ZZ} B_{in}.
$$

\begin{thm}
[Backelin-Fr\"oberg]
\label{thm.BF}
\cite[Prop. 3]{BF}
If $A$ is a connected graded $k$-algebra with defining relations of degree $\le d+1$, then $A^{(n)}$
is a quadratic algebra for all $n \ge d$.
\end{thm}

\begin{thm}
[Verevkin]
\label{thm.V}
 \cite[Thm. 4.4]{V}
Let $A$ be a connected graded algebra generated by $A_1$. Then $\QGr A \equiv \QGr A^{(n)}$ 
for all positive integers $n$.
\end{thm}

\begin{prop}
If Theorem \ref{thm.main0} holds for connected graded quadratic monomial algebras it holds for all connected graded monomial algebras.
\end{prop}
\begin{pf}
Let $A$ be a monomial algebra and give $\ell$, $F$ and $G$ the meanings they have in section \ref{sect.Ufn}.

By Theorem \ref{thm.BF}, $A^{(\ell)}$ is a quadratic algebra. Because $A$ is a monomial algebra so is 
$A^{(\ell)}$. By Theorem \ref{thm.V},   $\QGr A \equiv \QGr A^{(\ell)}$. Hence if Theorem \ref{thm.main0}
holds for $A^{(\ell)}$, then $\QGr A \equiv \QGr kQ'$ where $Q'$ is the Ufnarovskii graph for $A^{(\ell)}$. 
\end{pf}

\subsection{}
\label{ssect.example}
The Ufnarovskii graph for $A^{(\ell)}$ is more complicated than that for $A$.  For example, the Ufnarovskii graph 
for  $A=k\langle x,y\rangle/(y^3)$ is
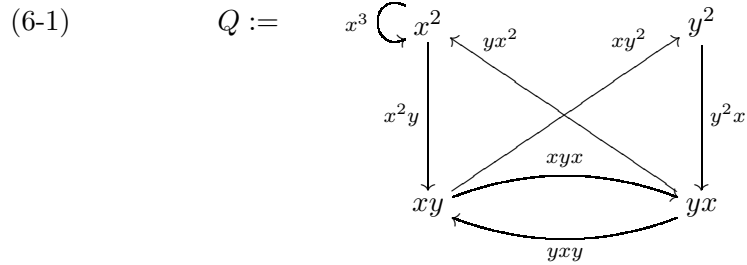
\begin{equation}
Q:= \phantom{xxx}  \UseComputerModernTips
\xymatrix{ 
\ar@(ul,dl)[]_{x^3}   x^2    \ar[dd]_{x^2y} &&& y^2  \ar[dd]^{y^2x} 
\\
\\
xy   \ar@/^1pc/[rrr]^{xyx} \ar[uurrr]^(0.8){xy^2}   &&& yx  \ar[uulll]_(0.8){yx^2} \ar@/^1pc/[lll]^{yxy}
}
\end{equation}
where the arrows are denoted by $w$ rather than $a_w$. The homomorphism $\fol:A \to kQ$ is given by
\begin{align*}
\fol(x)&=a_{x^3}+a_{x^2y}+a_{xyx} +a_{xy^2}
\\
\fol(y)&=a_{yx^2}+a_{yxy}+a_{y^2x}
\end{align*}
The $2$-Veronese subalgebra of $A$ is generated by  $s=x^2$, $t=xy$, $u=yx$, and $v=y^2$. We have
$$
A^{(2)}\cong \frac{k\langle s,t,u,v\rangle}{(vu,tv,v^2)}
$$ 
so its Ufnarovskii graph is 
$$
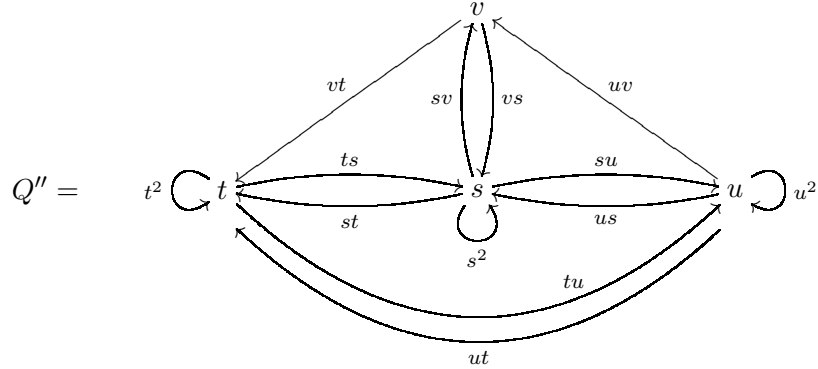
\UseComputerModernTips
\xymatrix{
&{} & {} & {} & {v} \ar[ddlll]_{vt} \ar@/^/[dd]^{vs} & {} & {} & {}
 \\
  \\
Q''= \phantom{xxx}
 &
{ t}\ar@(ul,dl)[]_{t^2} \ar@/_4pc/[rrrrrr]^(0.7){tu} \ar@/^/[rrr]^{ts} & {} & {} & {s}\ar@(dl,dr)[]_{s^2} \ar@/^/[lll]^{st} \ar@/^/[rrr]^{su} \ar@/^/[uu]^{sv} & {} & {} & {u}\ar@(ur,dr)[]^{u^2} \ar[uulll]_{uv} \ar@<2ex>@/^4pc/[llllll]^{ut} \ar@/^/[lll]^{us} }
$$
The homomorphism $f:k\langle s,t,u,v\rangle/(vu,tv,v^2) \to kQ''$ is given by
\begin{eqnarray*}
\fol(s)&=&a_{s^2}+a_{st}+a_{su}+a_{sv} \\
\fol(t)&=&a_{t^2}+a_{ts}+a_{tu} \\
\fol(u)&=&a_{u^2}+a_{us}+a_{ut}+a_{uv} \\
\fol(v)&=&a_{vs}+a_{vt}. 
\end{eqnarray*}
\section{A Remark}

The results in \cite{Sm2} and \cite{Sm3} show that many different $Q$ give rise to the equivalent categories 
$\QGr kQ$. 
Thus, given a finitely presented connected graded monomial algebra $A$, the Ufnarovskii graph is not the only $Q$ for which $\QGr A$ is equivalent to $\QGr kQ$. 

Consider, in particular, 
$$
A=\frac{k\langle x,y\rangle}{(y^3)}
$$
The Ufnarovskii graphs for $A$ and $A^{(2)}$ appear in section \ref{ssect.example}.
Since $A^{(\ell)}$ is quadratic for all $\ell \ge 2$, $\QGr kQ(A)\equiv \QGr kQ(A^{(\ell)})$  for all $\ell \ge 2$.

Furthermore, by \cite{Sm0}, $\QGr A$ is also equivalent to $\QGr kQ'$ where  
\begin{equation}
\label{Penrose.quiver}
Q'= \phantom{xxx}
  \UseComputerModernTips
\xymatrix{
0 \ar@(ul,dl)[]  \ar@/^.4pc/[r] & 1 \ar@/^.5pc/[l] \ar@/^.4pc/[r] & 2   \ar@/^1.3pc/[ll]
}
\end{equation}
\vskip .1in

\noindent
There is a direct proof of the equivalence $\QGr kQ(A)\equiv \QGr kQ'$.

\begin{thm}
 \cite{Sm3}
\label{thm1}
Let $L$ and $R$ be $\NN$-valued matrices such that $LR$ and $RL$ make sense.
Let $Q^{LR}$ be the quiver with incidence matrix $LR$ and $Q^{RL}$ the quiver with incidence matrix $RL$.
There is an equivalence of categories
$$
\QGr kQ^{LR} \equiv \QGr kQ^{RL}.
$$
\end{thm}

The equivalence $\QGr kQ(A)\equiv \QGr kQ'$ follows from Theorem \ref{thm1} because
$Q(A)=Q^{LR}$ and $Q'=Q^{RL}$ where
$$
L=\begin{pmatrix} 
1 & 0 & 0 \\
1 & 0 & 0 \\
0 & 1 & 0 \\
0 & 1 & 1 
\end{pmatrix}
\qquad \hbox{and} \qquad
R=\begin{pmatrix} 
1 & 0 & 0 & 1 \\
0 & 1 & 0 & 0 \\
0 & 0 & 1 & 0  
\end{pmatrix}.
$$


\end{document}